\documentclass[english,ams]{jcm}

\usepackage{babel}
\usepackage[latin1]{inputenc}

\usepackage{algorithm}
\usepackage{algorithmic}

\newtheorem{Def}{Definition}[section]

\newtheorem{Theo}[Def]{Theorem}
\newtheorem{Prop}[Def]{Proposition}
\newtheorem{Lem}[Def]{Lemma}
\newtheorem{Cor}[Def]{Corollary}
\theoremstyle{definition}
\newtheorem{remark}[Def]{Remark}

\newcommand{\F}{\mathbb{F}}
\newcommand{\N}{\mathbb{N}}

\newcommand{\ord}{\mathrm{ord}}
\newcommand{\ann}{\mathrm{ann}}
\newcommand{\len}{\mathrm{length}}
\newcommand{\lc}{\mathrm{lc}}
\newcommand{\GL}{\mathrm{GL}}
\newcommand{\lcm}{\mathrm{lcm}}

\newcommand{\Prob}{\mathrm{Prob}}
\newcommand{\calY}{\mathcal{Y}}
\newcommand{\rsp}{\mathrm{RowSp}}
\newcommand{\la}{\left<}
\newcommand{\ra}{\right>}
\newcommand{\ve}{\varepsilon}

\newcommand{\proofbeg}{\noindent\emph{Proof:}\ }
\newcommand{\proofof}[1]{\noindent\emph{Proof of #1:}}
\newcommand{\proofend}{\hfill$\blacksquare$}

\begin{document}

\title{Computing Minimal Polynomials of Matrices}

\author{Max Neunh\"offer}
\email{neunhoef@mcs.st-and.ac.uk}
\address{University of St Andrews, School of Mathematics and
Statistics,\\
North Haugh, St Andrews, Fife KY16 9SS, Scotland, UK}
\author{Cheryl E.~Praeger}
\email{praeger@maths.uwa.edu.au}
\address{University of Western Australia, School of Mathematics and
Statistics\\(M019), 35 Stirling Highway, Crawley 6009, Western Australia, Australia}

\begin{abstract}
We present and analyse a Monte-Carlo algorithm to compute the minimal 
polynomial of an $n\times n$ matrix over a finite field that requires
$O(n^3)$ field operations and $O(n)$ random vectors, and is well suited for
successful practical implementation. The algorithm, 
and its complexity analysis, use standard algorithms for polynomial 
and matrix operations. We compare features of the algorithm with 
several other algorithms in the literature. 
In addition we present a deterministic verification procedure
which is similarly efficient in most cases but has a worst-case
complexity of $O(n^4)$. Finally, we report the results of practical experiments
with an implementation of our algorithms in comparison with the
current algorithms in the {\sf GAP} library.
\end{abstract}

\classification{Primary: 15A21; Secondary: 15A15}

\keywords{Minimal polynomial, Frobenius normal form, Monte-Carlo algorithm,
randomisation, matrix, finite field}


\maketitle

\section{Introduction}

Let $\F$ be a finite field and $M \in \F^{n \times n}$ a matrix. This paper
presents and analyses a Monte Carlo  algorithm to compute the minimal 
polynomial of $M$, that is,
the monic polynomial $\mu \in \F[x]$ of least degree, such that
$\mu(M) = 0$. 
Determining the minimal polynomial is one of the fundamental computational
problems for matrices and has a wide range of applications. As well as 
revealing information about the Frobenius
normal form of $M$, the minimal polynomial also elucidates the structure 
of $\F^n$ viewed as $\F[x]$-module, where $x$ acts by multiplication with $M$. 
In addition the order of $M$ modulo scalars is often found by first 
determining the minimal polynomial. Apart from these applications it has
important practical utility, for example
in the context of the matrix group recognition project~\cite{OB}.

For these and other reasons a number of algorithms to determine the 
minimal polynomial may be found in the literature. We discuss some 
of them below. Our primary objective 
was to provide a simple and practical algorithm that could be 
implemented easily and would work well over small finite fields. In
particular we did not want to produce matrices with entries in 
larger fields or polynomial rings as intermediate results, 
and we preferred to restrict ourselves to using only row operations 
(rather than a combination of row and column operations). 
In addition we wished to use standard field and polynomial arithmetic,
and we wished to give  an explicit worst-case upper bound
for the number of elementary field operations needed, and not only an
asymptotic complexity statement.
Our Monte Carlo algorithm adheres to these requirements for matrices 
over fields $\F_q$ of order $q$.

\begin{Theo}\label{main}
For a given matrix $M \in \F_q^{n \times n}$ and a positive 
real number $\varepsilon < 1/2$, Algorithm~\ref{algminpolymc}
computes the minimal polynomial of $M$ with probability at least $1-\varepsilon$.
For sufficiently large $n$ and fixed $\varepsilon$, the number of elementary 
field operations required is less than $7n^3$ plus the costs of 
factorising a degree $n$ polynomial over $\F_q$ and constructing at most $n$
random vectors in $\F^n$.
\end{Theo}


Our algorithm to compute the minimal polynomial first computes the
characteristic polynomial in a standard way by spinning up and then
factoring out cyclic
subspaces. However, the novel aspect in this first phase
is the introduction of randomisation. While
not necessary for the computation of the characteristic polynomial
it underpins our proof of the Monte Carlo 
nature of our minimal polynomial algorithm. In addition to the Monte Carlo
minimal polynomial algorithm
we present and analyse in Section~\ref{verify} a deterministic verification procedure to
be run after Algorithm~\ref{algminpolymc} that has a similar asymptotic
complexity in many cases, but is $O(n^4)$ in the worst-case scenario.
Our motivation for giving concrete upper bounds for the costs of 
various component procedures was that, in practical implementations, 
these assist us to compare different algorithms in order to 
decide which to use in different situations. At the end of the paper we discuss
a practical implementation and tests of the algorithms in the {\sf GAP}
system~\cite{GAP4}.

\subsection{Other algorithms in the light of our requirements}

There are several interesting and asymptotically efficient
minimal polynomial algorithms for $n\times n$ matrices in the literature. The most  
asymptotically efficient deterministic algorithm is due to Storjohann 
\cite{Stor01} in 2001. It is nearly optimal, `requiring about the same number of 
field operations as required for matrix 
multiplication' (see \cite[Abstract, p368]{Stor01}). 
It involves a divide-and-conquer strategy that produces matrices 
with entries in polynomial rings as intermediate results. 
Changing the scalars to a larger field or polynomial ring is 
something we wished to avoid as it creates additional complications in practical 
applications within a computer algebra system used for group and 
matrix algebra computations.

Storjohann's earlier deterministic algorithm~\cite{Stor98} in 1998 uses classical 
field arithmetic and requires $O(n^3)$ field operations. 
It first reduces the matrix to `zig-zag form', using a mix of row 
and column operations, then produces the Smith normal form 
as a matrix with polynomial entries, and finally the Frobenius normal form.
In systems such as {\sf GAP}, matrices over small finite fields 
are stored in a compressed form that
makes row operations simple, but column operations difficult.
Restricting to one of these types of operations was one of our criteria.

A Monte Carlo minimal polynomial algorithm of Giesbrecht~\cite{Gie95} from 1995
that runs in `nearly optimal time' contains some 
features we find desirable for practical implementation, namely 
his algorithm first constructs a `modular cyclic decomposition' 
using random vectors, similar to our characteristic polynomial 
computation in Section~\ref{charpoly}. However, further steps include a 
modification of the `divide-and-conquer' 
Keller-Gehrig algorithm \cite{KelG85} and lead to a 
Las Vegas algorithm that computes a Frobenius form over an extension 
field and then the minimal polynomial.
The field size over which the given matrix is written is assumed 
to be greater than $n^2$, and if this is not the case it is 
suggested that an embedding into a larger extension field be used.
Several of these features were undesirable for us.

In \cite[Section 4]{AC97} Augot and Camion propose a deterministic algorithm
to compute the minimal polynomial of a matrix which is to some extent
similar to our algorithm. It is deterministic with
complexity $O(n^3 + m^2 \cdot n^2)$ field operations, 
where $m$ is the number of blocks
in the shift Hessenberg form. They prove that the complexity is $O(n^3)$ in the
average case. However, in the worst case it is $O(n^4)$, and 
no constants are provided in the complexity estimates. Although the
principal approach of their algorithm is similar to ours, the details
differ very much from our algorithm and analysis.

An interesting commentary on various algorithms, together with some new 
algorithms is given by Eberly~\cite{Eb00}. Eberly (see Theorem 4.2
in \cite{Eb00}) gives in particular a randomised algorithm for matrices 
over small fields that produces output from which (amongst other things)
the minimal polynomial can be computed, at a cost of $O(n^3)$.
The papers \cite{Eb00,Gie95,Steel,Stor98,Stor01} contain references to other 
minimal polynomial algorithms.  In all of the algorithms mentioned 
the asymptotic complexity statements give no information about 
the constants involved. 

On a practical note,  the minimal polynomial algorithm implemented in the GAP library
is the one in \cite{Steel} and (although we have been unable to 
confirm this) we assume that this is the algorithm implemented 
in Magma \cite{Magma}. 

\subsection{Outline of the paper}
In Section~\ref{notation} we introduce our notation, in
Section~\ref{complexity} we cite a few complexity bounds for basic
algorithms. The next Section~\ref{ordpoly} introduces order polynomials
and derives a few results about them. Then we turn to the computation
of the characteristic polynomial in Section~\ref{charpoly}, since this
is the first step in our minimal polynomial algorithm, which is described
and analysed in Section~\ref{minpoly}. We explain and modify the 
well-known algorithm to compute characteristic polynomials by introducing
some randomisation, because this is later needed in the analysis of our
main Monte Carlo algorithm. In Section~\ref{probest} we give some
probability estimates that are also used later in the analysis.
The second last Section~\ref{verify}
covers the deterministic verification of the results of our
Monte Carlo algorithm. We describe in detail cases in which this
verification is efficient and when it has a worse complexity.
Finally, in Section~\ref{performance} we report on the performance
of an implementation of our algorithm, including runtimes in 
realistic applications.  We compare these
times with the current implementation for minimal polynomial computations
in the {\sf GAP} library (see \cite{GAP4}), and as mentioned above, we 
believe that Magma and {\sf GAP} are both using the algorithm in 
\cite{Steel}. We show that our algorithm performs
much better in important cases, and that our bounds on the computing
cost are reflected in practical experiments.

\section{Notation}
\label{notation}

Throughout the paper $\F$ will be a fixed field. Although we envisage
$\F$ to be a finite field for our applications, this is not necessary
for most of our results. However in the later sections we use some probability 
estimates from Section~\ref{probest} that are only valid for 
finite fields.

By an \emph{elementary field operation} 
we mean addition, subtraction, multiplication or division of two field 
elements.
In all our runtime bounds we will assume that one elementary
field operation takes a fixed amount of time and we simply count
the number of such operations occurring in our algorithms.

We denote the set of $(m \times n)$-matrices over $\F$ by $\F^{m \times n}$
and the set of row vectors of length $m$ by $\F^m$. For a vector
$v \in \F^m$ we write $v_i$ for its $i$-th component and for a matrix 
$M \in \F^{m \times n}$ we denote its $i$-th row, which is
a row vector of length $n$, by $M[i]$. We use ``row vector
times matrix'' operations, and in general right modules throughout.
If $V$ is a vector space over $\F$ and $W$ is a subspace, the
quotient space is denoted by $V/W$ and its cosets by
$v+W$ for $v \in V$. The $\F$-linear span of the vectors
$v^{(1)}, \ldots, v^{(k)} \in V$ is denoted by 
$\left< v^{(1)}, \ldots, v^{(k)}\right>_\F$.

If $M \in \F^{n \times n}$ is a matrix and $V = \F^n$, we have a
natural action of $M$ as an endomorphism of $V$ by right multiplication.
The same holds for every $M$-invariant subspace $W < V$ and for
the corresponding quotient space $V/W$. We describe such a situation
by saying that ``the matrix $M$ induces an action on the $\F$-vector space''
$V, W, V/W$ respectively.

Throughout, $\F[x]$ denotes
the polynomial ring over $\F$ in an indeterminate $x$. For a square matrix $M$ and
a polynomial $p \in \F[x]$ we denote the evaluation of $p$ at $M$
by $p(M)$.

Whenever a matrix $M$ induces an action on a vector space $U$, we
will view $U$ as a right $\F[x]$-module by letting $x$ act like $M$,
that is $v \cdot x := v\cdot M$ in the above examples. We denote the
characteristic polynomial of this action by $\chi_{M,U}$. That is,
$\chi_{M,U}$ is the characteristic polynomial 
of the $(\dim_\F(U) \times \dim_\F(U))$-matrix given by choosing 
a basis of $U$ and writing the
action of $M$ induced on $U$ as a matrix with respect to that basis.
We use the same convention analogously for the corresponding minimal
polynomial $\mu_{M,U}$. Furthermore, we denote the $\F[x]$-submodule of
$U$ generated by the vectors $u^{(1)}, \ldots, u^{(n)}$ by $\left< u^{(1)}, \ldots,
u^{(n)} \right>_M$.

We use the two functions 
\begin{equation}\label{si}
s^{(1)}(a,b) := \sum_{i=b+1}^a i\quad \mbox{and}\quad
s^{(2)}(a,b) := \sum_{i=b+1}^a i^2
\end{equation}
for complexity expressions.
Note that for $a > b > c$ we have $s^{(j)}(a,c) = s^{(j)}(a,b) + s^{(j)}(b,c)$
for $j \in {1,2}$ and 
\begin{eqnarray}
\label{formels1}
s^{(1)}(n,0) &=& s^{(1)}(n,-1) = \frac{n(n+1)}{2}
\qquad\mbox{and} \\
\label{formels2}
s^{(2)}(n,0) &=& s^{(2)}(n,-1) = \frac{n(n+1)(2n+1)}{6}.
\end{eqnarray}

For later complexity estimates we note the following inequalities.

\begin{Lem}[Some upper bounds]
\label{estimates}
If $n = \sum_{i=1}^k d_i$ for some $d_i \in \N \setminus\{0\}$ and
$s_j := \sum_{i=1}^j d_i$ we have
\[ \sum_{j=1}^k s_j \le \frac{n(n+1)}{2} \quad \mbox{and} \quad
   \sum_{j=1}^k s_j(s_j+1) 
   \le \frac{n(n+1)(n+2)}{3}. \]
\end{Lem}
\proofbeg
We claim that for fixed $n$ both expressions are maximal if and only if
all $d_i$ are equal to one. We leave it to the reader to check that
both totals increase if we replace $d_j$ in some sequence
$(d_i)_{1 \le i \le k}$ by the two numbers $a$ and $d_j-a$ resulting
in the new sequence $(d'_1, \ldots, d'_{k+1}) := 
(d_1, d_2, \ldots, d_{j-1}, a, d_j-a, d_{j+1}, \ldots, d_k)$ 
of length $k+1$.
%
\proofend

\section{Complexity bounds for basic algorithms}
\label{complexity}

In some algorithms presented in later sections we use greatest common divisors
of univariate polynomials. To analyse these algorithms we use the 
following bounds which arise from standard polynomial computation. 
We take this approach because the standard algorithms for polynomials
are good enough for our complexity estimates in applications and we do not 
need the asymptotically best algorithms, discussion of which may be found 
conveniently in \cite{vzG}. 



\begin{Prop}[Complexity of standard greatest common divisor algorithm]
\label{standardgcd}\mbox{}

Let $f,g \in \F[x]$ with $n := \deg f \ge \deg g =: m$, and $f = qg + r$
with $q,r \in \F[x]$ such that $r=0$ or $\deg r < \deg g$. Then there is an
algorithm to compute $q$ and $r$ that needs less than $2(m+1)(n-m+1)$ 
elementary field operations. 

Furthermore, there is an
algorithm to compute $\gcd(f,g)$ that needs less than
$2(m+1)(n+1)$ elementary field operations.
\end{Prop}

\begin{remark}{We intentionally give bounds here which are not best possible,
since we want the bound for the $\gcd$ computation to be symmetrical in
$m$ and $n$.}
\end{remark}

\proofof{Proposition \ref{standardgcd}}
Use polynomial division and the standard {\sc Gcd} algorithm
and count. See \cite[Section 2.4 and Section 3.3]{vzG} for
smaller bounds that imply our symmetric bounds.
%
\proofend

\subsection{Polynomial factorisation}\label{polyfactn}

Some of our algorithms return partially factorised polynomials
which facilitate later factorisation into
irreducible factors. However, since the extent of this partial factorisation
is difficult to estimate, we use in our
analyses the complexity of finding the complete 
factorisation of a polynomial over a finite field as a product of 
irreducibles. We need such factorisations in our main algorithm.
In keeping with our other methods we make use of standard
polynomial factorisation procedures.

Details can be found in Knuth~\cite[4.6.2]{knuth} of a deterministic 
polynomial factorisation algorithm inspired by an idea of Berlekamp. 
Its cost is polynomial in both the degree $n$ and field size 
$|\F|=q$, as it requires $O(q)$ computations
of greatest common divisors. Thus it works well only for $q$ small. 
There is available a randomised (Las Vegas) version of the procedure
which (for arbitrary $q$) will always
return accurately the number $r$ of irreducible factors 
of $f(x)\in\F[x]$, but for which there is a small non-zero probability that it will 
fail to find all the irreducible factors.
It involves the procedure {\sc RandomVector}, which is discussed further in Subsection~\ref{random}, to 
produce independent uniformly distributed random elements of an 
$n$-dimensional vector space over $\F$  for 
which a basis is known. Throughout the paper logarithms are always taken to base $2$.

\begin{remark}[\sc PolynomialFactorisation]\label{rem:polyfactn}
Suppose $f(x)\in \F[x]$ of degree $n\geq 1$ with $r$ irreducible 
factors (counting multiplicities) and, if $q$ is large, suppose that 
we are given a real number $\ve$ such that $0<\ve<1/2$. The  number 
fact$(n,q)$ of elementary field operations required to find a complete
set of irreducible factors of $f(x)$ is at most

\smallskip
\begin{tabular}{lp{1.8in}}
$8n^3 + (3qr+17\log q)n^2$& deterministic algorithm \\
$O\big((\log \,\ve^{-1})(\log n)(\xi_n $ $+ n^2\log^3 q) +n^3\log^2q\big)$& 
Las Vegas algorithm
\end{tabular}\hspace*{-2mm}

\noindent where $\xi_n$ is an upper bound for the cost of one run of\/
{\sc RandomVector} on $\F^n$. The Las Vegas algorithm may fail, but with probability
less than $\ve$.
\end{remark}

\section{Order polynomials}
\label{ordpoly}

Let $M$ be a matrix in $\F^{n \times n}$ that induces an action
on an $\F$-vector space $V$.

We briefly recall the definition of the term ``order polynomial'':

\begin{Def}[Order polynomial $\ord_M(v)$ and relative order polynomial]
{\rm
The \emph{order polynomial} $\ord_M(v)$ of a vector $v \in V$ is the
monic polynomial $p \in \F[x]$ of smallest degree such that $v \cdot
p(M) = 0 \in V$. In particular $\ord_M(0)=1$.

For an $M$-invariant subspace  $W < V$, the \emph{relative order polynomial}
$\ord_M(v+W)$ (of $v$ relative to $W$) is the order polynomial of the element
$v+W \in V/W$ with respect to the induced action of $M$ on $V/W$.
}
\end{Def}

\begin{remark}
If we consider $V$ as an $\F[x]$-module as in 
Section~\ref{notation}, then $p$ is the monic generator
of the annihilator $\ann_{\F[x]}(v)$ of $v$ in $\F[x]$.
\end{remark}

The following observation follows immediately from the definition above.

\begin{Lem}[Relative order polynomials]
\label{relorderpol}
For an $M$-invariant subspace  $W < V$ and $v\in V$, 
$\ord_M(v+W)$ is the monic polynomial $p \in \F[x]$ of smallest degree such 
that $v \cdot p(M) \in W$.
\end{Lem}

We now turn to the question of how one computes the order polynomial of
a vector $v \in V$. The basic idea is to apply the matrix
$M$ to the vector repeatedly computing a sequence $v, vM, vM^2, \ldots, vM^d$
until $vM^d$ is a linear combination
\[ 
vM^d = \sum_{i=0}^{d-1} a_i vM^i, 
\]
with $a_i \in \F$, for $0 \le i < d$. If $d$ is minimal such that
this is possible, we have
\[ \ord_M(v) = x^d - \sum_{i=0}^{d-1} a_i x^i. 
\]
Although this procedure is simple and well-known, we 
present it in order to make explicit the number of elementary field 
operations needed.
To this end we describe in detail the computation of solutions for
the systems of linear equations involved. 

\begin{Def}[Row semi echelon form]
    \label{rowsemiechelon}
{\rm
    A non-zero matrix $S = (S_{i,j}) \in \F^{m \times n}$ is in 
    \emph{row semi echelon form} if there are positive 
    integers $r \le m$ and $j_1,\ldots,j_r \le n$ such that,
    for each $i \le r$, $S_{i,j_i} = 1$ and $S_{k,j_i} = 0$ for all $k > i$, 
    and also $S_{k,j} = 0$ whenever $k > r$. For $i \le r$, column 
    $j_i$ is called the \emph{leading column of row $i$}, and we write
    $\lc(i) = j_i$. A sequence of vectors $u^{(1)},\ldots,u^{(m)} \in F^n$ is 
    said to be in semi echelon form if the matrix with rows 
    $u^{(1)},\ldots,u^{(m)}$ is in row semi echelon form.

}
\end{Def}

\noindent
Note that in Definition~\ref{rowsemiechelon} we do not assume 
$j_1 < j_2 < \cdots < j_r$ which is the usual condition for an echelon
form.

\begin{Def}[Semi echelon data sequence]
{\rm    Let $Y \in \F^{m \times n}$ be a matrix with $m \le n$ and of rank $m$. A 
    \emph{semi echelon data sequence for $Y$} is a tuple $\calY=
(Y,S,T,l)$, where
    $S \in \F^{m \times n}$ is in row semi echelon form with
    leading column indices $l = (\lc(1), \ldots, \lc(m))$, and
    $T \in \GL(m,\F)$ with $TY=S$. Further, $T$ is a lower triangular
    matrix, that is, for $T = (T_{i,j})$ we have
    $T_{i,j} = 0$ for $i < j$. For a semi echelon data sequence $\calY$
    we call the number $m$ its \emph{length}, sometimes denoted
    $\len(\calY)$.
A semi echelon data sequence $\calY'=(Y',S',T',l')$ is said to 
\emph{extend} $\calY$
if $\len(\calY')>\len(\calY)$, the first $\len(\calY)$ rows of $Y'$ and $S'$ form
the matrices $Y$ and $S$ respectively, and the first $\len(\calY)$ entries of $l'$
form the sequence $l$. 
}
\end{Def}

\begin{remark}\label{rem:seds}
(a) The idea of this concept is that for a matrix 
$S \in \F^{m\times n}$ in row semi echelon form it is relatively
cheap to decide whether a given vector $v \in \F^n$ lies in the row
space of $S$, and if so, to write it as a linear combination of the
rows of $S$, that is, to find a vector $a \in \F^m$ such that
$v = aS = aTY$ (see Algorithm~\ref{clean}). Thus, the vector
$v$ is expressed as a linear combination of the rows of $Y$
using the vector $aT$ as coefficients.

(b) We call a semi echelon data sequence \emph{trivial} if $m=0$.
In this case, by convention, we take the row spaces of
the empty matrices $Y$ and $S$ to be the zero subspace of $\F^n$, we
denote the empty sequence in $\F^0$ by $0$, and we interpret $aS$ as
the zero vector of $\F^n$.
\end{remark}

\medskip
We now present Algorithm~\ref{clean}, which is one step in the
computation of a semi echelon data sequence for a matrix $Y$. We
denote by $S[i]$ the $i$-th row of the matrix $S$, and by $\rsp(S)$ the 
row space of $S$.

\begin{algorithm}
\caption{$\quad$ \sc CleanAndExtend}
\label{clean}
\begin{algorithmic}
\STATE \textbf{Input:} A semi echelon data sequence $\calY=(Y,S,T,l)$ with 
         $Y,S \in \F^{m \times n}$, $v \in \F^n$ (possibly $m=0$).
\STATE \textbf{Output:} A triple $(c,\calY',a')$ where $c$ is {\sc True} 
if $v \in \rsp(Y)$ and {\sc False} otherwise, $\calY'$ equals or 
extends $\calY$ respectively with $\len(Y') \le \len(Y)+1$, 
and $a'\in\F^{\len (\calY')}$, such that $v = a'S'$.
\vspace*{2mm}
\STATE $w := v$
\STATE $a := 0 \in \F^m$
\COMMENT{note that $w=v-aS$}
\FOR {$i = 1$ to $m$}
    \STATE $a_i := w_{l_i}$
    \STATE $w := w - a_i \cdot S[i]$
\ENDFOR\hspace*{6mm}
\COMMENT{still $w=v-aS$} 
\IF {w = 0}
    \STATE \textbf{return} $(\mbox{\sc True},(Y,S,T,l),a)$
\ELSE
    \STATE $j := $ index of first non-zero entry in $w$
    \STATE $a' := [ a \ \ w_j ]$,\quad $l' := \left[ l \ \ j \right]$, $\quad$
    \STATE $Y' := \left[ 
       \begin{array}{c} Y \\ v \end{array} \right]$, $\quad$
           $S' := \left[ 
       \begin{array}{c} S \\ w_j^{-1} \cdot w \end{array} \right]$, $\quad$
           $T' := \left[ 
       \begin{array}{cc} T & 0 \\ -w_j^{-1}\cdot aT & w_j^{-1} 
       \end{array} \right]$
    \STATE \textbf{return} $(\mbox{\sc False}, (Y',S',T',l'), a' )$
\ENDIF
\end{algorithmic}
\end{algorithm}

\begin{Prop}[Correctness and complexity of Alg.~\ref{clean}:
\label{PropCleanAndExtend}{\sc CleanAndExtend}]
The output of Algorithm~\ref{clean} satisfies the Output
specifications.
Moreover, Algorithm~\ref{clean} requires at most $2mn$ field operations
if $v\in\rsp(Y)$, and
$(2m+1)n + (m+1)^2 + 1$ field operations otherwise.
\end{Prop}

\begin{remark}\label{rem:cae}
(a) Given a semi echelon data sequence $(Y,S,T,l)$ with $Y,S \in \F^{m
\times n}$ and a vector $v \in \F^n$, Algorithm~\ref{clean} tries to
write $v$ as a linear combination of the rows of $S$. If this is not
possible, it constructs an extended semi echelon data sequence.

\medskip\noindent
(b) For the case of finite fields a simple and useful
optimisation is to reduce, where possible, the number of operations for vectors and matrices, 
for example, where a vector is multiplied by the zero scalar and the
result is added to some other vector. This can reduce the number of operations for
sparse vectors and matrices. Our estimates for the numbers of field operations then become
over-estimates.
\end{remark}

\smallskip
\proofof{Proposition~\ref{PropCleanAndExtend}}
The proof of the correctness of Algorithm~\ref{clean} is left to the
reader.
The \textbf{for} loop needs $2mn$ field operations if we count both multiplications
and additions. If $v\in\rsp(Y)$ then the algorithm terminates after this loop.
On the other hand, if  $v\not\in\rsp(Y)$, then 
Algorithm~\ref{clean} needs one inversion
of the scalar $w_j$ plus $2 \cdot \sum_{i=1}^m i = m(m+1)$
field operations for the vector times matrix multiplication
$aT$, because $T$ is a lower triangular matrix. This is altogether $m(m+1)+1$ operations. 
Finally, the scalar negation of $w_j^{-1}$ and the multiplication of $aT$ by 
$-w_j^{-1}$ needs another $m+1$ field operations, and a further $n$ operations 
are needed for the computation of $w_j^{-1} w$ in $S'$. Thus the
total number of field operations is at most $2mn + (m+1)^2 + 1 + n$.
\proofend

\smallskip
Having Algorithm~\ref{clean} at hand we can now present
Algorithm~\ref{algordpoly}, which computes relative order polynomials.
Since a (non-relative) order polynomial may be regarded as a relative order
polynomial with respect to the zero subspace, Algorithm~\ref{algordpoly}
can also be used to compute order polynomials, starting with the trivial semi
echelon data sequence, (see Remark~\ref{rem:seds}~(b)).

\begin{algorithm}[t]
\caption{$\quad$ \sc RelativeOrdPoly}
\label{algordpoly}
\begin{algorithmic}
\STATE \textbf{Input:} A semi echelon data sequence $\calY=(Y,S,T,l)$ with 
$Y,S \in \F^{m \times n}$ (possibly $m=0$),
$v \in \F^n$, 
and $M \in \F^{n \times n}$ such that $W:=\rsp(Y)$ is $M$-invariant.
\STATE \textbf{Output:} A triple $(p,\calY',b)$ consisting of the
relative order polynomial $p := \ord_M(v+W)$ of degree $d$,
a semi echelon data sequence $\calY'=(Y',S',T',l')$ of length $m+d$ equal
to or extending $\calY$, and a vector $b \in \F^{m+d}$ such that $vM^d =
bY'$.

\vspace*{2mm}
\STATE $(Y',S',T',l') := (Y,S,T,l)$ \hspace*{2mm}
\COMMENT{the primed variables change during the algorithm}
\STATE $v' := v$
\STATE $m' := m$ \hspace*{3.08cm}  \COMMENT{can be zero!}
\LOOP
    \STATE $(c,(Y',S',T',l'),a) := \mbox{\sc CleanAndExtend}((Y',S',T',l'),v')$  
    \STATE \hspace*{5cm} 
\COMMENT{$T'Y'=S'$, $v'=aS'$}
    \IF { $ c = \mbox{\sc True} $ }
        \STATE \textbf{leave loop}
    \ENDIF
    \STATE $v' := v' \cdot M$
     \STATE $m' := m' +1$	
\ENDLOOP \hspace*{1.5cm}  \COMMENT{at this stage $c=$ {\sc True}, $v'=aS'$, $m'=\len(\calY')$}
\STATE $d := m'-m$
\STATE $b := a\cdot T'$
\STATE $p :=x^d-\sum_{i=0}^{d-1} b_{m+1+i} x^i$
\STATE \textbf{return} $(p, (Y',S',T',l'),b)$
\end{algorithmic}
\end{algorithm}

\begin{Prop}[Correctness and complexity of Alg.~\ref{algordpoly}:
{\sc RelativeOrdPoly}]
\label{proprelorderpol}

\mbox{}\par
Let $\calY=(Y,S,T,l)$ be a semi echelon data sequence with $Y,S \in \F^{m \times
n}$ (possibly $m=0$), $v \in \F^n$, and $M \in \F^{n \times n}$
such that $W:=\rsp(Y)$ is $M$-invariant. 
The output of Algorithm 2 satisfies the Output specifications, and
moreover if $d>0$, then
rows $m+1, \ldots, m+d$ of $Y'$ are equal to $v,vM,\ldots,vM^{d-1}$
respectively.
Algorithm~\ref{algordpoly} requires at most
\begin{eqnarray*}
2dn^2 &+& (n+2)d +2(m+d)n 
+ 2(n+1)s^{(1)}(m+d-1,m-1) +  \\
 &+& s^{(2)}(m+d-1,m-1) 
+ 2s^{(1)}(m+d,0)
\end{eqnarray*}
elementary field operations where $s^{(1)}$ and $s^{(2)}$ are the
functions defined in (\ref{si}).
\end{Prop}

\begin{remark}
Note that, if $d=0$ then $S'=S$, so $v=b'Y\in W$, and in this case $p=1$. 
Algorithm~\ref{algordpoly} successively considers the vectors $v+W,
vM+W, \ldots, vM^d+W$  (those are the successive values of
$v'$) until $vM^d+W$ lies in the subspace of $V/W$ 
generated by the vectors $v+W, vM+W, \ldots, vM^{d-1}+W$. 
The given matrix $S$ together with Algorithm~\ref{clean} defines a
direct sum decomposition of the $\F$-vector space $V := \F^n = W \oplus W'$
where $W'$ is the subspace of vectors having $0$ in all positions
occurring in the list $l$. Since $W' \cong V/W$, 
Algorithm~\ref{algordpoly} effectively computes in $V/W$ by always `cleaning
out' vectors using $S$ first. 
\end{remark}

\proofof{Proposition~\ref{proprelorderpol}} 
We again leave the proof of correctness of Algorithm~\ref{algordpoly}
to the reader.
Algorithm~\ref{algordpoly} calls Algorithm~\ref{clean} ({\sc CleanAnd\-Extend}) exactly
$d+1$ times with the lengths of the input semi echelon data
sequences being $m,m+1, \ldots, m+d$. After each but the last call to
Algorithm~\ref{clean} the value of 
$c$ returned is {\sc False}, and after the last call the value of $c$ 
is {\sc True}. Thus, by
Proposition~\ref{PropCleanAndExtend}, the number of steps needed 
for the $d+1$ runs of Algorithm~\ref{clean} is at most
\[
\left(\sum_{i=m}^{m+d-1} ((2i+1)n+(i+1)^2+1)\right)  +  2(m+d)n  \]
\[
   = s^{(2)}(m+d-1,m-1) + (2n+2)s^{(1)}(m+d-1,m-1) + (n+2)d +2(m+d)n.
\]
In addition, we have to do $d$ multiplications of $v'$ with $M$, which require
$2n^2$ elementary field operations each, and finally the computation of $b$
requires $2s^{(1)}(m+d,0)$ elementary field operations, again since 
$T'$ is a lower triangular matrix. Summing up gives the expression in 
the statement. 
\proofend

\smallskip
We conclude this section with two lemmas that are used to
compute absolute order polynomials using relative ones. We again
view $V$ as an $\F[x]$-module by letting $x$ act like $M$. For
$\{v^{(1)},\dots,v^{(m)}\} \subseteq V$, we denote by $\la
v^{(1)},\dots,v^{(m)}\ra_M$ the submodule of $V$ generated
by $\{v^{(1)},\dots,v^{(m)}\}$, that is, the smallest $M$-invariant
subspace containing $\{v^{(1)},\dots,v^{(m)}\}$. If $m=1$ then $\la
v^{(1)}\ra_M$ is the $\F$-span of the set $\{v^{(1)},v^{(1)}M,\dots,
v^{(1)}M^{n-1}\}$. We call $\la v^{(1)}\ra_M$ a \emph{cyclic
subspace} relative to $M$.

\begin{Lem}[Order polynomials in cyclic subspaces]
\label{ordpolcyclic}
Let $v\in V$, $W = \left< v \right>_M < V$, and
$p := \ord_M(v)$ with $d := \deg(p)$. 
Then for each $w \in W$, there is a unique polynomial
$q \in \F[x]$ of degree less than $d$ such that $w = vq(M)$.
Moreover,
\[ \ord_M(w) = \frac{p}{\gcd(p,q)}. \]
\end{Lem}
\noindent We omit the routine proof for the sake of brevity. \proofend
%
%

%

\begin{Lem}[Absolute and relative order polynomials]
\label{absordpoly}
Let $W$ be an $M$-invariant subspace of $V$, $v \in V$ and 
$q := \ord_M(v+W) \in \F[x]$. Then
\[ \ord_M(v) = q \cdot \ord_M(vq(M)). \]
\end{Lem}
\noindent We omit the routine proof for the sake of brevity. \proofend

\section{Computing the characteristic polynomial}

\label{charpoly}

In this section we present a version of a standard algorithm for computing
the characteristic polynomial of a matrix together with its analysis. 
It differs from the standard version in its use of randomisation. 

\subsection{Random vectors}\label{random}

Our characteristic polynomial algorithm, and later ones, 
make use of the algorithms {\sc RandomVector} and 
{\sc RandomVector*} that produce independent uniformly 
dis\-tri\-bu\-ted random vectors, and independent uniformly distributed
random non-zero vectors,             
respectively, in a given finite vector space for which a basis is known.    
The algorithms are invoked for spaces
$\F^s$, for $s \in \N$, and for subspaces of $V$ of the form
\[
V(l)=\{ v\,|\, v_{l_i}=0\quad\mbox{for}\quad  1\leq i\leq m\}
\quad\mbox{where}\quad l=(l_1,\dots,l_m).
\]
If $l$ is the empty sequence then $V(l)=V$.
For a semi echelon data sequence $\calY = (Y,S,T,l)$, the 
vector space $V$ is the sum $V=V(l) \oplus \rsp(S)$. 

If $b=\mbox{\sc RandomVector}(\F^{\len(\calY)})$, then $bS$ is a
uniformly distributed random vector of $\rsp(S)$. 
Moreover we assume that for the disjoint spaces 
$\F^{\len(\calY)}$ and $V(l)$ the 
algorithms {\sc RandomVector} and {\sc RandomVector*} are applied 
independently so that in particular, 
if $a=\mbox{\sc RandomVector*}(V(l))$ 
then the sum 
$a + bS$ is a uniformly distributed random vector of $V\setminus \rsp(S)$. 

{\sc RandomVector} and, if we neglect the possibility of obtaining
the zero vector, also {\sc RandomVector*}, could proceed by selecting
independent uniformly distributed random field elements as coefficients
of the basis vectors. For the subspace $V(l)$, we could put zeros into
the entries occurring in $l$ and make random selections of elements from
$\F$ for each entry not in $l$.
For this reason we denote by $\xi_{r}$ an upper bound 
for the cost of {\sc RandomVector} or {\sc RandomVector*} applied to an 
$r$-dimensional space for one of these cases.
If $r<s$ then $\xi_r\leq\xi_s$ and $\xi_{r_1}+\xi_{r_2} \leq \xi_{r_1+r_2}$, and we would 
expect $\xi_r$ to vary linearly with $r$. In 
practical implementations the cost is much less than the cost of the field operations involved 
in the algorithm below.

\subsection{Characteristic polynomial algorithm}

\noindent The characteristic polynomial algorithm below would terminate           
successfully without making random selections of vectors. However, the  
use of randomisation is key to our application of this algorithm for    
finding minimal polynomials. As in previous sections, let $M$ be a      
matrix in $F^{n \times n}$ acting naturally on $V := \F^n$.

\begin{algorithm}
\caption{$\quad$ \sc CharPoly}
\label{algcharpoly}
\begin{algorithmic}
\STATE \textbf{Input:} $M \in \F^{n\times n}$
\STATE \textbf{Output:} 
A tuple $(k, (p^{(j)})_{1\leq j\leq k}, \calY, (b^{(j)})_{1\leq j\leq
k})$, where  each $p^{(j)}\in\F[x]$ and $\prod_{i=1}^k p^{(i)} =
\chi_{M,V}$ is the characteristic polynomial
of $M$ in its action on $V$,  each  $b^{(j)}\in\F^n$ and $\calY$ is a
semi echelon data sequence of length $n$ with the properties
specified in Proposition~\ref{propcharpoly}.

\vspace*{2mm}
\STATE $i := 0$
\STATE $\calY^{(0)}:=$ a trivial semi echelon data sequence
\WHILE {$\len(\calY^{(i)}) < n$}
    \STATE $i := i + 1$
    \STATE $a := \mbox{\sc RandomVector}(\F^{\len(\calY^{(i-1)})})$
    \STATE $c := \mbox{\sc RandomVector*}(V(l^{(i-1)}))$
    \STATE $v^{(i)} := aS^{(i-1)} + c$
	\STATE 
       \hspace*{10mm} \COMMENT{ $v^{(i)}\not\in\rsp(S^{(i-1)})$ where $\calY^{(i-1)}=(Y^{(i-1)},S^{(i-1)},
                              T^{(i-1)},l^{(i-1)})$}
    \STATE $(p^{(i)},\calY^{(i)},b^{(i)}) :=$
             $\mbox{\sc RelativeOrdPoly}(\calY^{(i-1)},v^{(i)},M)$
\STATE  \hspace*{10mm} \COMMENT{$b^{(i)}\in\F^{\len(\calY^{(i)})}$; we add $n-\len(\calY^{(i)})$ zeros to make $b^{(i)}\in\F^{n}$}
\ENDWHILE
\STATE $k := i$
\STATE \textbf{return} $(k,(p^{(j)})_{1 \le j \le k}, 
                       \calY^{(k)}, 
(b^{(j)})_{1 \le j \le k})$
\end{algorithmic}
\end{algorithm}

\begin{Prop}[Correctness and complexity of Algorithm~\ref{algcharpoly}]
\label{propcharpoly}\mbox{}
Algorithm~\ref{algcharpoly} satisfies
the Output specifications, and furthermore
$\calY=(Y,S,T,l)$ where
$Y \in \F^{n \times n}$ is invertible with rows
\[ v\!^{(1)}, v\!^{(1)}M, \ldots, v\!^{(1)}M^{d_1-1}, 
v\!^{(2)}, v\!^{(2)}M, \ldots, v\!^{(2)}M^{d_2-1},
   \ldots, v\!^{(k)}, v\!^{(k)}M, \ldots, v\!^{(k)} M^{d_k-1} 
\]
where $d_i := \deg(p^{(i)})$ for $1 \le i \le k$. Further, for $1 \le i \le k$,
$v^{(i)}$ is a uniformly distributed random element of\/ $V\setminus W_{i-1}$, 
$v^{(i)} M^{d_i} = b^{(i)} Y$ and $p^{(i)} = \ord_M(v^{(i)} + W_{i-1})$,
where $W_{i-1} := \left< v^{(1)}, \ldots, v^{(i-1)}\right>_M$ (for $i>1$), 
an $M$-invariant subspace of $V$ of dimension $s_{i-1}:=\sum_{j=1}^{i-1} d_j$, and 
$W_0=0$ of dimension $s_0=0$.
Moreover, Algorithm~\ref{algcharpoly} requires at most
\[
\frac{33}{6}n^3+4n^2+\frac{3}{2}n 
\]
elementary field operations, plus
$k\xi_{n}$ for the $k$ calls to {\sc RandomVector*} and {\sc RandomVector}. Neglecting 
the latter cost this  is less than $6n^3$ elementary field operations, for sufficiently large $n$.
\end{Prop}

\begin{remark}
We denote the semi echelon data sequences 
$\calY^{(i)}$ in the algorithm using indices to enable us to speak 
more easily about the intermediate results. However in practice
we have only one variable $\calY=(Y,S,T,l)$, the entries of which are 
growing during the execution of the algorithm.
\end{remark}


\begin{remark}
Note that we do not multiply together the factors of 
$\chi_{M,V}$ because in our application of Algorithm~\ref{algcharpoly} 
we do not need the product itself.
\end{remark}

\proofof{Proposition~\ref{propcharpoly}}
Most statements in the proposition follow immediately from 
Proposition~\ref{proprelorderpol}: note that, because of the conventions
explained in Remark \ref{rem:seds}(b), in the first run of the
`while' loop $v^{(1)}=c$ is a random non-zero vector of $V$ and
$p^{(1)}=\ord)M(v^{(1)})$, and more generally,
in the $i^{th}$ run of the 
`while' loop, Algorithm~\ref{algcharpoly}
chooses a  vector $v^{(i)}$ that is a uniformly distributed random element 
of $V\setminus \rsp(S^{(i-1)})$
and applies Algorithm~\ref{algordpoly}.
This immediately establishes all statements about
$(Y,S,T,l)$ including the one about the invertibility and the
rows of $Y$. Also it is clear that $p^{(i)} = \ord_M(v^{(i)} + W_{i-1})$.

Next we show that $\prod_{i=1}^k p^{(i)} = \chi_{M,V}$. 
This follows by considering
the matrix $YMY^{-1}$, which has the same characteristic polynomial
as $M$. Considering the action of $M$ with respect to the 
ordered basis of $\F^n$ given by the rows of $Y$, it follows from
the construction that $YMY^{-1}$ 
(written with respect to the standard basis) is equal to
\[ \left[\begin{array}{cccc}
 C_1       &   0  & \cdots & 0 \\
 B^{(2)}_1 &  C_2 & \ddots & \vdots \\
 \vdots    &\ddots& \ddots & 0 \\
 B^{(k)}_1 &\cdots& B^{(k)}_{k-1} & C_k
\end{array} \right] \]
where the matrix $C_i$ is the companion matrix of the polynomial $p^{(i)}$,
and the $B^{(i)}_j$, for $2 \le i \le k$ and $1 \le j \le i-1$, are matrices
in $\F^{d_i \times d_j}$
with one non-zero row at the bottom and all other rows zero.
If $b^{(i)} = (b^{(i)}_1,\dots,b^{(i)}_n)$, then the bottom row of 
$B^{(i)}_j$ is $(b^{(i)}_{s_{j-1}+1},\dots,b^{(i)}_{s_j})$.
With this format at hand it is clear that the characteristic polynomial
of $YMY^{-1}$ is equal to the product $\prod_{i=1}^k p^{(i)}$ because
the $C_i$ are companion matrices.

Finally we derive the statement about the number of elementary field operations
needed by Algorithm~\ref{algcharpoly}.
In the $i^{th}$ run of the {\bf while} loop, the cost of constructing the
random vectors $a$ and $c$ is at most 
\[ \xi_{n-\len(\calY^{(i-1)})}+\xi_{\len(\calY^{(i-1)})}
\leq \xi_{n}, 
\]
(see Subsection~\ref{random}). The cost to compute $v^{(i)}$ is 
at most $2s_{i-1}n$ elementary field operations,
where $s_0=0$, and for $i\geq1$, $s_{i}=\sum_{j=1}^{i}d_j$ with 
$d_j=\deg p^{(j)}$. The cost 
of applying Algorithm {\sc RelativeOrdPoly}
is, by Proposition~\ref{proprelorderpol}, at most
\[ 
2d_in^2 + (n+2)d_i +2s_{i}n + 2(n+1)s^{(1)}(s_i-1,s_{i-1}-1) 
+ s^{(2)}(s_i-1,s_{i-1}-1) + 2s^{(1)}(s_i,0) 
\]
elementary field operations, noting that the value 
of `$d$' is $d_i$, the value of `$m$' is $s_{i-1}$, $s_{i-1}+d_i=s_i$,
and $s^{(1)}$, $s^{(2)}$ are the functions defined in (\ref{si}).

We consider the different terms
one by one, summing each over $i$ from $1$ to $k$. 
The total cost of constructing the random vectors is at most $k\xi_{n}$.
Summing the terms $2 s_{i-1} n$ gives $2n\sum_{i=1}^{k} s_{i-1}$, and
summing the terms $2d_in^2$ gives $2n^3$ since $\sum_{i=1}^k d_i=n$. 
Similarly, summing the terms
$(n+2)d_i$ gives $(n+2)n$. From the terms $2s_in$ we
get a contribution of 
$2n \sum_{i=1}^{k} s_i$.
The next two expressions involving the functions $s^{(1)}$ and $s^{(2)}$
sum to $2(n+1)s^{(1)}(n-1,0) = n(n+1)(n-1)$ and $s^{(2)}(n-1,0) = 
\frac{(n-1)n(2n-1)}{6}$ respectively, using (\ref{formels1}) and
(\ref{formels2}) and the properties noted
above it. 
Finally, the terms $2s^{(1)}(s_i,0)$ sum to 
$2\sum_{i=1}^k s^{(1)}(s_i,0) 
= \sum_{i=1}^k s_i(s_i+1)$.
Thus in total we obtain $k\xi_n$ plus
\begin{eqnarray*}
 2n^3
   &+&n(n+1)(n-1)
   +\frac{(n-1)n(2n-1)}{6}
   +n(n+2)
   +2n\sum_{i=1}^k s_{i-1} \\
   &+&2n\sum_{i=1}^k s_i
   +\sum_{i=1}^k s_i(s_i+1) 
\end{eqnarray*}
elementary field operations. The first four of these terms sum to $\frac{10}{3}n^3
   +\frac{3}{2}n^2
   +\frac{7}{6}n$. 
Using Lemma~\ref{estimates},
\[ 
   2n\sum_{i=1}^k s_{i-1}
   +2n\sum_{i=1}^k s_i
   +\sum_{i=1}^k s_i(s_i+1)\leq 2n^2(n+1)+\frac{n(n+1)(n+2)}{6} \]
so the total cost is at most 
\[
\frac{33}{6}n^3+4n^2+\frac{3}{2}n +k\xi_n.
\]
For sufficiently large $n$ this is less than $6n^3+k\xi_n$.
\proofend

\section{Probability estimates using the structure theory for modules}
\label{probest}

The basic idea of our minimal polynomial Algorithm~\ref{algminpolymc}
is to compute the order polynomials of a few
random vectors under the action of a given matrix $M$ and to prove that, 
with high probability, their least common multiple is
equal to the minimal polynomial of $M$.
The purpose of this section is to use the structure theory 
of $V = \F^n$ as an $\F[x]$-module to derive probability
estimates to be used in that proof.

First suppose that the characteristic polynomial of $M$ is written as a product
$\chi_{M,V} = \prod_{i=1}^t q_i^{e_i}$ with pairwise distinct
irreducible polynomials $q_i \in \F[x]$ and positive integer
multiplicities $e_i$.

Using \cite[Theorem 3.12]{Jacob1} we can then write the $\F[x]$-module $V$
as a direct sum of primary cyclic modules
\begin{equation} \label{primary}
V \cong \bigoplus_{i=1}^t \bigoplus_{j=1}^{m_i} w_{i,j} \F[x] 
\end{equation}
such that $\ord_M(w_{i,j}) = q_i^{f_{i,j}}$ with
$e_i \ge f_{i,1} \ge f_{i,2} \ge \cdots \ge f_{i,m_i} \ge 1$
and $\sum_{j=1}^{m_i} f_{i,j}=e_i$ for $1 \le i \le t$.

The minimal polynomial $\mu_{M,V}$ is the least
common multiple of the order polynomials of the vectors 
$(w_{i,j})_{1 \le i \le t, 1 \le j \le m_i}$, and hence is 
$\mu_{M,V} = \prod_{i=1}^t (q_i)^{f_{i,1}}$.

We use this structural description to derive the
first probability bound for the case where $\F = \F_q$ is a finite field
with $q$ elements.

\begin{Prop}[Probability that a 
$q_i$ has equal mult.~in $\mu_{M,V}$
and $\ord_M(v)$]
\label{ProbOneMult}

\mbox{}\par
Let $\F = \F_q$ be a finite field with $q$ elements, let $V=\F^n$, let
$U$ be a (possibly zero) $M$-invariant subspace such that the multiplicity of
$q_i$ in $\mu_{M,U}$ is strictly smaller than in $\mu_{M,V}$, and let 
$v$ be a uniformly distributed random element of $V\setminus U$. Then 
the multiplicity of $q_i$ is the same in $\ord_M(v)$ and $\mu_{M,V}$ 
with probability greater than
$1-q^{-\deg q_i}$.
\end{Prop}
\proofbeg
By assumption the multiplicity of $q_i$ in $\mu_{M,U}$ 
is less that its multiplicity $f:=f_{i,1}$ in $\mu_{M,V}$. 
Let $w:=w_{i,1}$, with $w_{i,1}$ as in (\ref{primary}), so that $V=X\oplus Y$ with
$X, Y$ invariant under $M$ and $X=\la w\ra_M$. Then $\mu_{M,X}=
q_i^{f}$. We may identify the primary cyclic $\F[x]$-module 
$X$ with  $w\F[x]$, which is
isomorphic to the module $\F[x]/(q_i^{f}\F[x])$,
and in turn this is uniserial with composition series
\[ 
0 <   \frac{q_i^{f-1}\F[x]}{q_i^{f}\F[x]}
     <   \frac{q_i^{f-2}\F[x]}{q_i^{f}\F[x]} <
\cdots <  \frac{q_i\F[x]}{q_i^{f}\F[x]} 
       < \frac{\F[x]}{q_i^{f}\F[x]}. 
\]
Thus, $X$ has a unique maximal $\F[x]$-submodule, namely  $X':=\la w q_i(M)\ra_M$,
and $X'$ has codimension $r:=\deg(q_i)$ in $X$. 

As discussed above, each vector $v\in V$ has a 
unique expression as $v=x+y$ with $x\in X, y\in Y$. Moreover
$\ord_M(v)$ is the least common multiple of $\ord_M(x)$ and $\ord_M(y)$.
In particular, if $x\not\in X'$, then $\ord_M(x)=q_i^{f}$ and hence
the multiplicity of $q_i$ in $\ord_M(v)$ and $\mu_{M,V}$ is the same.
The number of vectors $v=x+y$ with $x\not\in X'$ is 
\[
|X\setminus X'|\cdot |Y|=(1-\frac{1}{q^r})|X|\cdot|Y|=(1-\frac{1}{q^r})q^n.
\]
Each of these vectors $v$ lies in $V\setminus U$ since the multiplicity of
$q_i$ in $\mu_{M,U}$ is less than $f$. Thus the probability, for
a uniformly distributed random $v\in V\setminus U$, that the multiplicity 
of $q_i$ in $\ord_M(v)$ and $\mu_{M,V}$ is the same is at least
\[
(1-\frac{1}{q^r})\frac{q^n}{|V\setminus U|} > 1-\frac{1}{q^r}.
\]
\proofend

\begin{remark}
If for some irreducible factor $q_i$ we have $m_i > 1$
and $f_{i,1} = f_{i,2}$, then the above probability is even higher,
because we can apply the above argument independently to two or more summands
$w_{i,1}\F[x]$ and $w_{i,2}\F[x]$.
\end{remark}

\smallskip
We now give a second probability bound which will be crucial in our
Monte Carlo algorithm to compute the minimal polynomial. In 
that algorithm we choose  a sequence of vectors 
$v^{(1)}, \dots v^{(u)}$ such that $v^{(1)}$ is a uniformly distributed
random element of $V\setminus\{0\}$, and for $i\geq2$ we choose
$v^{(i)}$ as a uniformly distributed random element of 
$V \setminus U$, where $U=\left< v^{(1)}, \ldots, v^{(i-1)} \right>_M$.
We hope to find $\mu_{M,V}$ as the least common multiple of the
orders of these vectors.

\begin{Prop}[Probability that an lcm of order polynomials equals
$\mu_{M,V}$]
\label{ProbAllMult}

\mbox{}\par
Let\/ $\F = \F_q$ be a finite field with $q$ elements.
Suppose a sequence of vectors $v^{(1)}, \ldots, v^{(u)} \in V$ is chosen
as follows: $v^{(1)}$ is a uniformly
distributed random element of $V\setminus\{0\}$, and for $i>1$,
$v^{(i)}$ is a  uniformly distributed random element of  
$V \setminus \left< v^{(1)}, \ldots, v^{(i-1)} \right>_M$. 
Let
\[ 
f := \lcm( \ord_M(v^{(1)}), \ord_M(v^{(2)}), \ldots, 
\ord_M(v^{(u)}) ). 
\]
Then the probability that  $f = \mu_{M,V}$  is greater than
\[ 1-\sum_{i=1}^t q^{-u\deg q_i}. \]
\end{Prop}
\proofbeg
Consider the random experiment described in the statement. We first
examine one irreducible factor $q_i$. Let $E_i$ denote the event
that the multiplicity of $q_i$ in $f$ is strictly smaller than 
the multiplicity $f_{i,1}$ of $q_i$ in $\mu_{M,V}$. Furthermore, for
$1 \le j \le u$, let $F_j$  be the event that the multiplicity of $q_i$ in
$\ord_M(v^{(j)})$ is strictly smaller than $f_{i,1}$. 

Note that
the $F_j$ are not stochastically independent since we choose
$v^{(j)}$ outside of $\left< v^{(1)}, \ldots, v^{(j-1)}\right>_M$.
However, $E_i = F_1 \cap F_2 \cap \cdots \cap F_u$ because $f$ is the
least common multiple of the order polynomials of the $v^{(j)}$.
By Proposition~\ref{ProbOneMult} applied with $U=\{0\}$, 
the probability $\Prob(F_1)$
is less than $q^{-\deg q_i}$. Moreover, in the situation that 
$F_1 \cap \cdots \cap F_j$ holds and $j<u$,
we apply Proposition~\ref{ProbOneMult} with the
subspace $U := \left< v^{(1)}, \ldots, v^{(j)} \right>_M$ to conclude
that the conditional probability $\Prob(F_{j+1} | F_1 \cap \cdots \cap F_j)$
is less than $q^{-\deg q_i}$.
Thus we have
\begin{eqnarray*}
\Prob(E_i) = \Prob(F_1)&\cdot& \Prob(F_2 | F_1) \cdot \Prob(F_3 | F_1 \cap F_2)
   \cdot \cdots \\
    && \cdots \cdot \Prob(F_u | F_1 \cap \cdots \cap F_{u-1}) 
   < q^{-u\deg q_i}.
\end{eqnarray*}
Finally we consider all the different irreducible factors $q_i$. 

Even though
the events $E_1,\ldots,E_t$ may not be stochastically independent, we have
\[ \Prob( E_1 \cup \cdots \cup E_t ) \le \sum_{i=1}^t \Prob(E_i)
   < \sum_{i=1}^t q^{-u\deg q_i} \]
as claimed.
\proofend

\section{Computing minimal polynomials}
\label{minpoly}

Our minimal polynomial algorithm runs Algorithm~\ref{algcharpoly} as its
first step. So assume, from now on,
that we have already run Algorithm~\ref{algcharpoly}
and obtained all the output it produces, in particular
the basis given by the rows of the matrix $Y$ (as in Proposition~\ref{propcharpoly}),
\[ 
(v^{(1)}, v^{(1)}M, \ldots, v^{(1)} M^{d_1-1}, \ldots, v^{(k)}, 
v^{(k)} M, \ldots, v^{(k)} M^{d_k-1}) 
\]
the relative order polynomials $p^{(i)} = \ord_M(v^{(i)}
+ W_{i-1})$, and the vectors $b^{(i)}$ for $1 \le i \le k$.
Also assume that we have factorised all the polynomials $p^{(i)}$
as products $p^{(i)} = \prod_{j=1}^t q_j^{e_{i,j}}$ of irreducible 
polynomials $(q_j)_{1 \le j \le t}$.

The matrices $M$ and $YMY^{-1}$ have the same characteristic
and minimal polynomials. Also the order polynomials $\ord_M(v)$ and
$\ord_{YMY^{-1}}(vY^{-1})$ are equal for every $v \in V$ and thus also
the order polynomials $\ord_M(vY)$ and $\ord_{YMY^{-1}}(v)$ are equal for every
$v \in V$.

For the convenience of the reader we display the matrix
$YMY^{-1}$ in Figure~\ref{bigmat}. Note in particular that the matrix is
sparse, provided that the degrees $d_i$ are not too small.
Due to the special form of $YMY^{-1}$ it is much more efficient to compute
the images of vectors under $YMY^{-1}$, than under $M$. 
This is crucial in the analysis of our algorithms. 
Therefore we will from now on do all computations of order polynomials
with respect to $YMY^{-1}$.

Set $M' := YMY^{-1}$ and $W'_i := W_i Y^{-1}$ for $1 \le i \le k$.
Note that we have $v^{(i)} = e^{(s_{i-1}+1)} Y$ for $1 \le i \le k$
where $e^{(1)}, \ldots, e^{(n)}$ is the standard basis of $\F^n$. (That is,  
 $e^{(i)}$ contains exactly one $1$ in position $i$ and otherwise zeros. Recall
that $s_i = \sum_{j=1}^i d_j$ with $s_0 = 0$.) Furthermore, for $1 \le i \le k$,
the space $W_i = \left< v^{(1)}, \ldots, v^{(i)} \right>_M$ is equal to the space
$\{ vY \mid v \in \F^n \mbox{ with } v_j = 0 \mbox{ for } j > s_i \}$. 
Thus, the space $W'_i$ is the $\F$-linear
span $\left< e^{(1)}, e^{(2)}, \ldots, e^{(s_i)}\right>_\F$ and we have a
filtration
\[ 0 = W'_0 < W'_1 < W'_2 < \cdots < W'_k = V \]
such that each quotient $W'_i/W'_{i-1}$ is an $M'$-cyclic space generated 
by the coset represented by the standard basis vector $e^{(s_{i-1}+1)}$.

We begin by presenting Algorithm~\ref{algordpolabs} which computes the
absolute order polynomial of a vector with respect to the matrix $YMY^{-1}$,
using all the data acquired during Algorithm~\ref{algcharpoly}. 
We will apply this later
in the minimal polynomial algorithm to the first few of the vectors 
$e^{(s_{i-1}+1)}$ 
produced during a run of Algorithm~\ref{algcharpoly}. Note that for 
the analysis it is crucial that a number $z$ such that
the vector $v$ lies in $W'_z$ 
is given as input to the algorithm.

\begin{figure}
\caption{Overview of the matrix $YMY^{-1}$}
\label{bigmat}
\[ \left[ \begin{array}{cccccc|cccccc|c|cccccc}
\cline{1-6}
  0 & 1 &      &      &   &   \\
    & 0 & 1    &      & 0 &   \\
    &   &\ddots&\ddots&   &   \\
    & 0 &      &   0  & 1 &   \\
    &   &      &      & 0 & 1 \\
  {}* & * &   *  &   *  & * & * \\
\cline{1-12}
  &&&&&& 0 & 1 &      &      &   &   \\
  &&&&&&  & 0 & 1    &      & 0 &   \\
  &&&&&&  &   &\ddots&\ddots&   &   \\
  &&&&&&  & 0 &      &   0  & 1 &   \\
  &&&&&&  &   &      &      & 0 & 1 \\
  {}*&*&*&*&*&*&* & * &   *  &   *  & * & * \\
\cline{1-13}
  &&&\vdots&&& &&&\vdots&&& \ddots \\
\cline{1-19}
  &&&&&& &&&&&& & 0 & 1 &      &      &   &   \\
  &&&&&& &&&&&& & & 0 & 1    &      & 0 &   \\
  &&&&&& &&&&&& & &   &\ddots&\ddots&   &   \\
  &&&&&& &&&&&& & & 0 &      &   0  & 1 &   \\
  &&&&&& &&&&&& & &   &      &      & 0 & 1 \\
  {}*&*&*&*&*&*& *&*&*&*&*&*&* &*& * &   *  &   *  & * & * \\
\hline
\end{array} \right] \]
\end{figure}

\begin{algorithm}
\caption{$\quad$ \sc OrdPoly}
\label{algordpolabs}
\begin{algorithmic}
\STATE \textbf{Input:} $M$, $k$, $(Y,S,T,l)$, $(p^{(j)})_{1 \le j \le k}$,
$(b^{(j)})_{1 \le j \le k}$ as returned by {\sc CharPoly}, an integer
$z$ with $1 \le z \le k$,
$v \in W'_z$, and the factorisation  
$p^{(j)} = \prod_{r=1}^t q_r^{e_{j,r}}$ for all $j\leq k$
\STATE \textbf{Output:} A list of factorised polynomials, the product
of which is $\ord_{YMY^{-1}}(v)$ 
\vspace*{2mm}
\STATE $i := z$ \hspace*{1cm} \COMMENT{will run down to $1$}
\STATE $f := [\ ]$ \hspace*{1cm} \COMMENT{empty list}
\REPEAT
    \STATE $h := \sum_{j=1}^{d_i} v_{s_{i-1}+j} x^{j-1}$
    \IF {$h \neq 0$}
    \STATE $\hat g := p^{(i)}/\gcd(h,p^{(i)})$ \hspace*{1cm} \COMMENT{factorised}
        \STATE \textbf{add} $\hat g$ to list $f$
        \STATE \textbf{compute} product $g$ of factors in $\hat g$
        \IF {$i > 1$}
            \STATE $v := v \cdot g(YMY^{-1})$ \hspace*{1cm} 
            \COMMENT{see Proposition~\ref{propordpol} for this computation}
        \ENDIF
    \ENDIF
    \STATE $i := i - 1$
\UNTIL {$i=0$}
\STATE \textbf{return} $f$
\end{algorithmic}
\end{algorithm}

\begin{Prop}[Correctness and complexity of Algorithm~\ref{algordpolabs}:
{\sc OrdPoly}]
\label{propordpol}

\mbox{}\par
Let $\F = \F_q$ be a field with $q$ elements.
The output of Algorithm~\ref{algordpolabs} satisfies the Output specifications.
Moreover, Algorithm~\ref{algordpolabs}
requires at most
\[
\sum_{j=1}^{z} 
  \left( 4 d_j^2  + 3d_j s_j + 2d_j \sum_{r=1}^j s_r \right)
\ \leq\ (\frac{z}{2}+9)s_z^2
\]
elementary field operations, 
where $d_j=\deg p^{(j)}$, $s_j=\sum_{r=1}^j d_j$ for $j\geq1$ and
$s_0=0$; and this is less than $n^3$ for $n$ sufficiently large.

\end{Prop}
\proofbeg
Since we are computing an order polynomial with respect to the matrix
$M' = YMY^{-1}$ we can always use the form of this matrix as displayed
in Figure~\ref{bigmat}.

The basic idea of Algorithm~\ref{algordpolabs} is to use
Lemmas~\ref{ordpolcyclic} and \ref{absordpoly} applied to the filtration
\[ 0 = W'_0 < W'_1 < W'_2 < \cdots < W'_k = V. \]

Starting with $i := z$ and the original $v$ lying in the space 
$W'_z$, the variable $i$ runs downwards until $1$.
In each step, Algorithm~\ref{algordpolabs}
computes the relative order polynomial $g := \ord_{M'}(v+W'_{i-1})$
for the then current value of $v \in W'_i$.
This assertion follows from
Lemma~\ref{ordpolcyclic} noting that,
by our discussion above, $p^{(i)} = \ord_M(v^{(i)}+W_i)= \ord_{M'}(e^{(s_{i-1}+1)}+W'_{i-1})$.
Next $v g(M')$ is evaluated, which lies in $W'_{i-1}$ by
Lemma~\ref{relorderpol}, and the induction can go on with $i$
replaced by $i-1$. The product of the polynomials in the list $f$ 
returned is the product of all 
the relative order polynomials computed in the repeat loop, and this is equal to
$\ord_{M'}(v)$, by Lemma~\ref{absordpoly}.

To count the number of elementary field operations is a bit complicated
here. Note first that by assumption we already know a factorisation of all the
$p^{(i)} = \prod_{j=1}^t q_j^{e_{i,j}}$ into
irreducible factors. Now  $\gcd(h,p^{(i)})$ is equal to the product of
the greatest common divisors $\gcd(h,q_j^{e_{i,j}})$, for $j\leq t$. Since the degrees of the
polynomials  $q_j^{e_{i,j}}$ sum up to the degree $d_i$ of $p^{(i)}$, finding these gcd's, by 
Proposition~\ref{standardgcd}, requires at most
\[ 
2(\deg(h)+1) \cdot \sum_{j=1}^t \left(\deg(q_j) e_{i,j} + 1\right) 
\] 
field operations, which is at most $4d_i^2$ since
$\deg(q_j) e_{i,j} + 1 \le 2 \deg(q_j) e_{i,j}$. Note that this is a rather
crude estimate. At this stage we know all multiplicities of the
$q_j$ in $\gcd(h,p^{(i)})$ and thus in $g := p^{(i)}/\gcd(h,p^{(i)})$.
Thus we have computed $g$ in factorised form, which is denoted by
$\hat g$ in Algorithm~\ref{algordpolabs}.

Now we discuss the number of operations needed to evaluate $v g(M')$.
Due to the sparseness of $M'$, a multiplication of a vector $w$ of
$W_i'$ from the right by $M'$ needs only a shift (which we neglect here)
and an addition of a multiple of the non-zero
part of $b^{(r)}$ for $1 \le r \le i$ requiring $2s_r$ operations for
each $r$. Thus computing $wM'$ requires at most $\sum_{r=1}^i 2s_r $
elementary operations. Note that $wM'$ lies in $W_i'$ still. If
$f(x)\in\F_q[x]$ of degree $d$, say $f(x)=\sum_{r=0}^dc_rx^r$, then we
can compute $wf(M')=\sum_{r=0}^dc_rw(M')^r$ by first computing
$w(M')^r\in W_i'$ for $1\leq r\leq d$ at a cost of at most
$2d\sum_{r=1}^i s_r $, next computing $c_rw(M')^r$ for $0\leq r\leq d$
at a cost of at most $(d+1)s_i$, and then adding these vectors at a
further cost of at most
$ds_i$, making a total cost to compute $wf(M')$ of at most
$(2d+1)s_i +2d\sum_{r=1}^i s_r $ elementary field operations.

The polynomial $g(x)$ is available in factorised form,
say $g(x)=\prod_{s=1}^uf_s(x)$ with $\deg f_s=m_s$,
where $\sum_{s=1}^um_s=\deg g=d_i$. From the previous
paragraph we see that $vg(M')$ can be computed at a cost of at most
\[
\sum_{s=1}^u\left((2m_s+1)s_i +2m_s\sum_{r=1}^i s_r \right)
=(2d_i+u)s_i+2d_i\sum_{r=1}^is_r \leq 3d_is_i+2d_i\sum_{r=1}^is_r.
\]
Subsequent runs of the repeat loop require similar numbers of elementary
operations, with $i$ replaced by $j$ where $i-1\geq j\geq 1$.
Thus Algorithm 4 needs at most
\[
\sum_{j=1}^z \left( 4d_j^2  + 3d_j s_j + 2d_j 
\sum_{r=1}^j s_r \right)
\]
elementary field operations, as claimed in the proposition.

To find a simpler upper bound we look at the terms one by one. The
last and most important term can be bounded by
\begin{eqnarray*}
  2\sum_{j=1}^z \left( d_j \sum_{r=1}^j s_r \right)
    &=& 2 \sum_{r=1}^z s_r \left(\sum_{j=r}^z d_j\right)
    = 2 \sum_{r=1}^z s_r (s_z - s_{r-1}) \\
    &=& 2 \sum_{r=1}^z s_r (s_z - s_r) + 2 \sum_{r=1}^z s_r d_r
    \le (\frac{z}{2}+2) \cdot s_z^2
\end{eqnarray*}
since the function $t(s_z-t)$ has maximum value $s_z^2/4$
for $t$ in the interval $[0,s_z]$.

The second term $3\sum_{j=1}^z d_j s_j$ is at most $3s_z^2$,
and the term
$4 \sum_{j=1}^z d_j^2$ is at most $4s_z\sum_{j=1}^z d_j = 4s_z^2$.

Altogether this amounts to a bound of $(\frac{z}{2}+9)s_z^2$ as claimed.
Asymptotically, this is bounded above by $n^3$ in the worst case as
$n \to \infty$.
\proofend

\smallskip
Now we present our main procedure, Algorithm~\ref{algminpolymc}. 

\begin{algorithm}
\caption{$\quad$ \sc MinPolyMC}
\label{algminpolymc}
\begin{algorithmic}
\STATE \textbf{Input:} $M \in \F_q^{n \times n}$, $\varepsilon$ with
$0 < \varepsilon < 1/2$.
\STATE \textbf{Output:} A tuple $(b,f)$ where $b$ is either {\sc True}
or {\sc Uncertain} and $f \in \F_q[x]$
\STATE \hspace*{0mm} \phantom{\textbf{Output:}} (see
Proposition~\ref{propminpoly} for details).
\vspace*{2mm}
\STATE
       $((p^{(j)})_{1 \le j \le k},(Y,S,T,l),(b^{(j)})_{1 \le j \le k})
       := \mbox{\sc CharPoly}(M)$
\STATE Factorise all $p^{(j)} = \prod_{r=1}^t q_r^{e_{j,r}}$ 
\STATE Determine the least $u \in \N$ such that
 $\sum_{r=1}^t q^{-u\deg q_r} \le \varepsilon$
\STATE $u := \min\{ u,k \}$
\STATE $f := \lcm(p^{(1)}, \ldots, p^{(k)})$ 
\FOR {$i = 2$ to $u$}
    \STATE $f := \lcm(f,\mbox{\sc OrdPoly}
           (M,k,(Y,S,T,l),(p^{(j)})_{1 \le j \le k}, 
           (b^{(j)})_{1 \le j \le k},i,e^{(s_{i-1}+1)}))$
\ENDFOR
\IF {$u = k$ or $\deg f = n$}
    \STATE \textbf{return} $(\mbox{\sc True},f)$
\ELSE
    \STATE \textbf{return} $(\mbox{\sc Uncertain},f)$
\ENDIF
\end{algorithmic}
\end{algorithm}

\begin{Prop}[Correctness and complexity of Algorithm~\ref{algminpolymc}:
{\sc MinPolyMC}]\label{propminpoly}

\mbox{}\par
Given a matrix $M \in \F_q^{n \times n}$ and a number
$\varepsilon$ with $0 < \varepsilon < 1/2$, Algorithm~\ref{algminpolymc}
returns a tuple $(b,f)$, where $b$ is either {\sc True} or {\sc Uncertain}
and $f \in \F_q[x]$ is a polynomial. With probability at least $1-\varepsilon$
the polynomial $f=\mu_{M,\F_q^n}$, and if
\/ $b = \mbox{\sc True}$ then
$f = \mu_{M,\F_q^n}$ is guaranteed. Moreover, 
if\/ $f \ne \mu_{M,\F_q^n}$, then 
$f$ is a proper divisor of $\mu_{M,\F_q^n}$ and
every irreducible factor of $\mu_{M,\F_q^n}$ divides $f$.

The number of elementary field operations needed by 
Algorithm~\ref{algminpolymc} is bounded above by
\[ 
{\rm char}(n,q) + {\rm fact}(n,q) + 
  \sum_{i=1}^u (\frac{i}{2}+9)s_i^2
\]
where ${\rm char}(n,q)$ is an upper bound for the number of elementary field operations needed to
compute the characteristic polynomial (see
Proposition~\ref{propcharpoly}), ${\rm fact}(n,q)$ is an upper bound for 
the number of elementary field
operations needed to factorise each of a set of polynomials over $\F_q$ whose degrees sum to $n$
(see Subsection~{\rm\ref{polyfactn}}). Moreover either $u=k$, or $u<k$ and
$\sum_{j=1}^t q^{-u\deg q_j} \le \varepsilon$.

For $n$ sufficiently large and fixed $\varepsilon$, this is less than
\[ 6n^3 + {\rm fact}(n,q) + 
\frac{1}{3} \lceil\frac{\log n-\log \varepsilon}{\log q}\rceil^2 \cdot n^2 
\]
which is less than $7n^3 + {\rm fact}(n,q)$ (plus the cost of computing at most $n$ random 
vectors in Algorithm~\ref{algcharpoly}).
\end{Prop}

\begin{remark}
Note that if we use a randomised polynomial factorisation algorithm
(necessary for large $q$), then the algorithm can be modified to allow
for a possible failure of factorisation of the `Factorise' step (line
2). Thus Theorem~\ref{main} follows from Proposition~\ref{algminpolymc}.
An upper bound for the term {\rm fact}$(n,q)$ in the complexity
bound is given in Remark~{\rm\ref{rem:polyfactn}}, and this yields an
upper bound in Proposition~{\rm\ref{propminpoly}} of $O(n^3\log^3 q)$
for $n$ sufficiently large and fixed $\varepsilon$.
\end{remark}

\begin{remark}
    \label{algdeterm}
    Algorithm~\ref{algminpolymc} can be changed into a deterministic
    one by running the `$i$-loop' for $i$ up to $k$, instead of $u$. An upper
    bound of the cost is then given by replacing $u$ by $k$ in the
    formula for the cost in Proposition~\ref{propminpoly}.

If $k>u$, then these additional $k-u$ runs of the `$i$-loop' may be viewed as 
a `verification  algorithm'. By Proposition~\ref{algminpolymc}, 
the additional cost of these extra runs is
\[
\sum_{i=u+1}^k(\frac{i}{2}+9)s_i^2 \leq s_k^2\left(\frac{(k+u+1)(k-u)}{4}+9(k-u)\right)
= s_k^2\frac{(k-u)(k+u+37)}{4}
\]
and for sufficiently large $n$ this cost is less than 
$n^4/4$ field operations.
\end{remark}

\proofof{Proposition~\ref{propminpoly}}
Algorithm~\ref{algminpolymc} first computes the characteristic
polynomial of $M$ in its action on $\F_q^n$ and its factorisation.
This computation provides firstly the irreducible factors $q_j$ of
the minimal polynomial that allow us to determine $u$, and secondly
the input needed for running Algorithm~\ref{algordpolabs} to compute
the order polynomials of $v^{(2)}, \ldots, 
v^{(u)}$. Thirdly, it also yields a nice base change matrix $Y$ such 
that these order polynomials with respect to the matrix $M$ can in
fact be determined using Algorithm~\ref{algordpoly} for the vectors
$e^{(s_1+1)},\ldots,e^{(s_{u-1}+1)}$ since we have
$\ord_M(v^{(i)}) = \ord_{YMY^{-1}}(e^{(s_{i-1}+1})$ for $1 \le i \le
k$. Note that
$p^{(1)} = \ord_{YMY^{-1}}(e^{(1)}) = \ord_M(v^{(1)})$. 
By Proposition~\ref{propcharpoly},
the vector $v^{(j)}$ is a uniformly distributed random element of
$V\setminus\{0\}$ if $j=1$, or $V\setminus\left< v^{(1)}, \ldots,
v^{(j-1)}\right>_M$ if $j>1$. Hence, by Propositions~\ref{ProbAllMult}
and~\ref{propordpol}, the probability that $f$ after termination of
Algorithm~\ref{algminpolymc} is equal to $\mu_{M,\F_q^n}$ is at least
$1-\varepsilon$.

  From the discussion at the beginning of Section~\ref{probest},
$\mu_{M,\F_q^n}$ is the least common multiple of the $k$ polynomials
$\ord_M(v^{(1)}), \dots,\ord_M(v^{(k)})$, and hence if $u=k$ then
$f = \mu_{M,\F_q^n}$. This also implies, since the initial value
of $f$ is $\lcm(p^{(1)}, \ldots, p^{(k)})$, that the returned
polynomial $f$ divides $\mu_{M,\F_q^n}$ and every irreducible factor of
$\mu_{M,\F_q^n}$ divides $f$. In particular, if $\deg f =n$ then we must
have $f=\chi_{M,\F_q^n}=\mu_{M,\F_q^n}$. Thus if $(\mbox{\sc True},f)$
is returned then $f = \mu_{M,\F_q^n}$ is guaranteed.

The number of elementary field operations needed follows from 
Propositions~\ref{propcharpoly} and~\ref{propordpol} and summing. Note that,
after the factorisations computed in line 2 of the algorithm, we neglect
the forming of least common multiples and the products here, because
all results from Algorithm~\ref{algordpolabs} come already factorised
into irreducible factors. We can thus compute the least common multiples
by taking maximums of multiplicities. Hence the first displayed upper
bound is proved.

For the asymptotic complexity bound we have to consider the initial value of the number $u$,
namely the least integer $u$ such that $\sum_{j=1}^t q^{-u \deg q_j} \le \varepsilon$.
The largest value of this sum occurs when all the $q_j$ have degree 1, and as there are 
then at most $n$ such polynomials,  $\sum_{j=1}^t q^{-u \deg q_j}\le nq^{-u}$. 
Thus $u$ is at most the least integer such that $nq^{-u}\le \varepsilon$, namely 
\[
u_0:=\lceil \frac{\log n + \log (\varepsilon^{-1})}{\log q}\rceil
\]
and the value of $u$ used in the algorithm is at most $\min\{u_0,k\}\leq u_0$.
By Proposition~\ref{propcharpoly}, the asymptotic value of char$(n,q)$ is
less than $6n^3$ for $n$ sufficiently large, (plus the cost $k\xi_n$ 
of making $k$ random selections of vectors). 
By Proposition~\ref{propordpol},  the number of elementary
field operations used for the computation of the  $u\le u_0$ order polynomials
is at most
\[
\sum_{i=1}^u(\frac{i}{2}+9)s_i^2 \leq \frac{u(u+1)}{4}s_u^2+9us_u^2
\le u_0\left(\frac{u_0+37}{4}\right) n^2.
\]
which, for sufficiently large $n$ and fixed $\epsilon$, is less than
$\frac{1}{3}u_0^2n^2<n^3$.
\proofend

\section{Deterministic verification}
\label{verify}

In this section we explain how the probabilistic result of our Monte Carlo
algorithm can be verified deterministically. We begin by discussing cases
that can be handled rather cheaply, before we present several general verification procedures, all of which, unfortunately, have a worst-case cost of $O(n^4)$ field operations.

All notation from previous sections remains in force. The first result 
follows immediately from Proposition~\ref{propminpoly}.

\begin{Prop}[Cases, in which the result is already proven to be correct]
If the output polynomial of Algorithm~\ref{algminpolymc} is  $\chi_{M,\F_q^n}$,
then the output is $(${\sc True}, $\chi_{M,\F_q^n})$
and is correct.
\end{Prop}

For the next result observe that, if Algorithm~\ref{algminpolymc} 
is modified so that the {\bf for} loop is run $k$ times, then the resulting 
polynomial $f$ is guaranteed to be the minimal polynomial, giving a deterministic
algorithm with proven result. (Proof of correctness 
is given in the proof of Proposition~\ref{propminpoly}.)

\begin{Prop}[Case of few random vectors chosen during comp.~of
$\chi_{M,\F_q^n}$]
\label{veryfewvectors}

\mbox{}\par
If $k \le \sqrt{n}$,  and the {\bf for} loop in 
Algorithm~\ref{algminpolymc} is run $k$ times, then the output polynomial
is $\mu_{M,\F_q^n}$. The 
overall cost of this modification of Algorithm~\ref{algminpolymc} is at most
\[ {\rm char}(n,q) + {\rm fact}(n,q) + \frac{1}{4}n^3+\frac{37}{4} n^{5/2} \]
elementary field operations.
\end{Prop}
\proofbeg The only change 
to the complexity estimate is for the number of elementary field operations 
in the second last line of the proof of Proposition~\ref{propminpoly}:
\[
\sum_{i=1}^k(\frac{i}{2}+9)s_i^2 \leq \frac{k(k+1)}{4}s_k^2+9ks_k^2
\leq \frac{k(k+1)}{4}n^2+9kn^2\le \frac{1}{4}n^3+\frac{37}{4} n^{5/2}.
\]
The rest follows from Proposition~\ref{propminpoly}.
\proofend

\medskip
For the case of large $k$, we may use the procedure suggested in Remark\ref{algdeterm} as a verification algorithm, at a cost of $O(k^2n^2)$ field operations. Two alternative verification procedures are given below. The first involving evaluation on vectors is given in Proposition~\ref{eval}, and the second using null space computations is given in Proposition~\ref{propverify}.

\begin{Prop}[Verification by evaluation on vectors]\label{eval}
    \mbox{}\par
    For the output $(\textsc{Uncertain},f)$ of Algorithm~\ref{algminpolymc} 
    one can verify $f = \mu_{M,\F_q^n}$ using
    at most $dn(k-u)(k+u+4)$ elementary field operations where $u$ and
    $k$ are as in Proposition~\ref{propminpoly} and $d = \deg f$.
\end{Prop}
\proofbeg
The idea here is to check whether $e^{(s_{i-1}+1)} f(YMY^{-1})$ is equal to zero, for 
$u+1 \le i \le k$,
by direct evaluation using the techniques described in
the proof of Propostion~\ref{propordpol}. Recall first that the
result $f$ comes in factorised form. Since $e^{(s_{i-1}+1)}$ lies
in $W'_i$ the arguments in the proof of Proposition~\ref{propordpol}
show that this evaluation can be done using at most
$3d s_i + 2d \sum_{r=1}^i s_r$ elementary field operations. Thus, an
upper bound for the total cost for all these evaluations is
\[ 3d \sum_{i=u+1}^k s_i + 2d \sum_{i=u+1}^k \sum_{r=1}^i s_r. \]
The first term is bounded above by $3dn(k-u)$. As to the second
term, for $1 \le j \le u+1$, the value $s_j$ occurs in this expression
with coefficient $2d(k-u)$, while for $u+2 \le j \le k$, it occurs
with coefficient $2d(k-j+1)$. Thus the second term is at most
\[ 
2dn(u+1)(k-u) + 2dn(k-u)(k-u-1)/2 = dn(k-u)(k+u+1). 
\]
Adding this to the upper bound for the first term we get at most
$dn(k-u)(k+u+4)$ as claimed.
\proofend

\medskip
For the following discussion we need a lemma:

\begin{Lem}[Cost of evaluation of a polynomial at a matrix]
\label{costpolyeval}
Let $M \in \F^{n \times n}$ be a matrix and $f \in \F[x]$ a polynomial 
with degree $d < n$. Then the evaluation $f(M)$ can be computed using
at most $2dn^3$ elementary field operations.
\end{Lem}

\proofbeg We take $2n^3$ elementary field operations as an upper bound for a matrix
multiplication. The computation of the powers $M^2, M^3, \ldots, M^d$
needs at most $2(d-1)n^3$ elementary field operations. The
multiplication, for each $i=1,\dots,d$, of $M^i$ by a coefficient of $f$
and addition of the result to the already computed matrix (the sum of previous terms) 
needs another $2dn^2$ elementary field operations. Finally, the
constant term of $f$ has to be added along the diagonal, which is yet
another $n$ elementary field operations. Since $d+1\le n \le 2 n^2$, 
this is altogether at most $2dn^3$
as claimed.
\proofend

\smallskip
Of course, this immediately implies:

\begin{Cor}[Small degree minimal polynomial]
If $\deg \mu_{M,\F_q^n} < n$, then the output of Algorithm~\ref{algminpolymc} 
can be verified by evaluation using at most 
\[ 2\cdot n^3 \cdot \deg \mu_{M,\F_q^n} \] 
elementary field operations.
\end{Cor}

\begin{remark}
Note that using \cite[Theorem 2]{AC97} we could lower the complexity
in Lemma~\ref{costpolyeval} to $O(\sqrt d n^3)$ provided we 
stored $O(\sqrt d)$ matrices in memory at the same time. However, since
storing a matrix in $\F^{n \times n}$ needs $O(n^2)$ of memory, this
approach would often become impractical before  a
concrete problem would become intractable because of time constraints. 
We use our estimates in Lemma~\ref{costpolyeval} because of
these practical considerations.
However, in some practical situations, an improved polynomial evaluation
algorithm using more memory may be suitable.
\end{remark}

We now present Algorithm~\ref{algminpolyverify} that can be run after
Algorithm~\ref{algminpolymc} to verify the correctness of the resulting polynomial
deterministically.

\begin{algorithm}
\caption{$\quad$ \sc MinPoly verification}
\label{algminpolyverify}
\begin{algorithmic}
\STATE \textbf{Input:} $M \in \F^{n \times n}$, $\chi_{M,V} = \prod_{i=1}^t q_i^{e_i}$ (factorised), 
and a candidate $\prod_{i=1}^t q_i^{f_i}$ for $\mu_{M,\F_q^n}$ (factorised),
all data from Algorithm~\ref{algminpolymc}
\STATE \textbf{Output:} {\sc True} or a positive number $j$ (see
Proposition~\ref{propverify} for details).
\vspace*{2mm}
\FOR {$i=1$ to $t$}
    \IF {$f_i < e_i$}
        \STATE $M' := q_i(YMY^{-1})^{f_i}$
        \STATE $d := \dim_\F( \ker(M') )$
        \IF {$d < \deg(q_i) \cdot e_i$} 
            \STATE \textbf{return} $i$
        \ENDIF
    \ENDIF
\ENDFOR
\STATE \textbf{return} {\sc True}
\end{algorithmic}
\end{algorithm}

\begin{Prop}[Deterministic minimal polynomial verification]
\label{propverify}\mbox{}

If Algorithm~\ref{algminpolyverify} is called with candidate minimal polynomial
$\prod_{i=1}^t q_i^{f_i}$ from Algorithm~\ref{algminpolymc}, 
then it either returns {\sc True} or a 
positive integer~$j$.
In the former case, $\mu_{M,\F_q^n}=\prod_{i=1}^t q_i^{f_i}$,
while in the latter case the
multiplicity of $q_j$ in  $\mu_{M,\F_q^n}$ is greater than $f_j$.
The number of elementary field operations required by 
Algorithm~\ref{algminpolyverify} is at most
\[ 
 n^3 \cdot \sum_{i=1}^t \left( 2\deg q_i + 2\lceil \log f_i \rceil+1 \right).
\] 
\end{Prop}

\proofbeg Let $r_i := \deg q_i$ for $i = 1, \ldots, t$.
We again view $\F^n$ as $\F[x]$-module as in Section~\ref{probest} by
letting $x$ act as right multiplication by $M$. By 
\cite[Theorem~3.12]{Jacob1}, it is isomorphic to a direct sum of 
primary cyclic $\F[x]$-modules
\[ 
\F^n \cong \bigoplus_{i=1}^t \bigoplus_{j=1}^{m_i} w_{i,j} \F[x], 
\]
such that $\ord_M(w_{i,j}) = q_i^{f_i,j}$ with 
$e_i \ge f_{i,1} \ge f_{i,2} \ge \cdots \ge f_{i,m_i} \ge 1$ and
$\sum_{j=1}^{m_i} f_{i,j} = e_i$. Thus, for each $i$,
$q_i$ occurs in  $\mu_{M,\F_q^n}$ with multiplicity $f_{i,1}$,
and so in particular $f_i\le f_{i,1}$.
The element $q_i^{f_i}$ acts
invertibly on all direct summands $w_{i',j} \F[x]$ with $i' \neq i$
since $q_i$ is irreducible and every order polynomial of a non-zero vector
in such a direct summand is a power of $q_{i'}$, by 
Lemma~\ref{ordpolcyclic}. For $i' = i$ however, the dimension of the kernel 
of the action of $q_i^{f_i}$ on $w_{i,j} \F[x]$ is 
$r_i \cdot \min\{f_i,f_{i,j}\}$.
Thus the dimension of the kernel of the action of $q_i^{f_i}$ on the
whole of $\F^n$ is equal to 
\[
r_i\sum_{j=1}^{m_i}\min\{f_i,f_{i,j}\}
\le r_i\sum_{j=1}^{m_i}f_{i,j}=r_ie_i
\]
with equality if and only if $f_i\ge f_{i,1}$. 
Since $f_i\le f_{i,1}$, equality holds above if and only if $f_i$
is equal to the multiplicity $f_{i,1}$ of $q_i$ in  $\mu_{M,\F_q^n}$.
Therefore, Algorithm~\ref{algminpolyverify} always returns the
result as stated in the Proposition.

As to the cost, Algorithm~\ref{algminpolyverify} evaluates $q_i$ at 
$YMY^{-1}$ which needs at most $2r_i n^3$ elementary field operations
by Lemma~\ref{costpolyeval}. It then takes the result to the
$f_i^{\mathrm{th}}$ power, which can be done by repeated squaring
with at most $2n^3\lceil \log f_i \rceil$ elementary field operations, and 
finally computes the dimension of a null space, which can be done with at most
$n^3$ elementary field operations (compute a semi echelon basis
of the row space of the matrix). Note that we are not using the 
sparseness of $YMY^{-1}$ here.
\proofend

\begin{remark}
The cost in Proposition~\ref{propverify} is much smaller than $n^4$ in
many cases. One of the worst cases occurs when $\chi_{M,\F^n}$ contains lots of
different factors of degree $1$ each occurring with multiplicity $3$, and 
all the $f_i$ are equal to $2$. Then Algorithm~\ref{algminpolyverify} has to
square about $n/3$ matrices and compute the null spaces of the results. 
This amounts to about
$2n^4/3$ elementary field operations, which is only about twice as fast as 
directly evaluating the minimal polynomial at $M$. Note that even in this case
only about every sixth entry of\/ $YMY^{-1}$ is different from zero. 
\end{remark}

\begin{remark}
 As in each of our procedures there are some simplifications we could make in practice
which do not reduce the worst case complexity estimates. For example, in Algorithm~\ref{algminpolyverify}, there is no need to compute the kernel of $q_i(YMY^{-1})^{f_i}$ if the irreducible $q_i$ does not divide any of the relative order polynomials $p^{(j)}$ for $u<j\leq k$. 
\end{remark}

\section{Performance in practice}
\label{performance}

In this section we give some experimental evidence concerning 
the performance of Algorithm~\ref{algminpolymc} in comparison with 
that of algorithms currently implemented in the
{\sf GAP} library (see \cite{GAP4}). 

All computations were done on a machine with an Intel Core 2 Quad CPU Q6600
running at 2.40 GHz with 8 GB of main memory and two times 4 MB of second level
cache. 

We were unable to confirm that Magma \cite{Magma} 
uses an algorithm  based on the canonical forms algorithm of Alan
Steel presented in~\cite{Steel} for computing minimal 
polynomials, although this is indicated in~\cite[Abstract]{Steel} and on the
web %
(see \texttt{http://magma.maths.usyd.edu.au/magma/htmlhelp/text347.htm}).

Our colleague Colva Roney-Dougal kindly
ran the Baby Monster example matrix $M_2$ on Magma and the resulting
times were roughly equivalent to the timing in the column ``Lib'' of
Figure~\ref{timings}, suggesting that this is indeed the case. Since the
minimal polynomial algorithm in the {\sf GAP} library is also based on 
the algorithm in \cite{Steel}, we did not conduct extensive comparison tests of our
algorithm on Magma.

\subsection{Guide to the test data}
The timing results
are in Figure~\ref{timings}, all times are in seconds. 
The column marked ``$n$'' contains the
dimension of the matrix, the column marked ``$q$'' the number of elements
of the base field. 
The columns marked ``Lib'' and ``AS'' contain the times needed for one run of 
the minimal polynomial algorithm based on \cite{Steel} as implemented in 
the {\sf GAP} library, and as implemented (by the first author) in the
{\sf GAP} language, respectively.
The column ``MC'' contains the total time for our 
Monte Carlo algorithm as presented in Algorithm~\ref{algminpolymc}.
The next three columns marked ``Spin'', ``Fact'' and ``OrdP''
contain the times for the three phases of this algorithm respectively, 
namely the first phase to compute the characteristic polynomial via
relative order polynomials, the second phase to factor all factors of
the characteristic polynomial and count multiplicities, and the third
phase to compute some absolute order polynomials to guess the minimal
polynomial. Finally, the last column marked ``Ver.'' contains
the time for the deterministic verification via 
Algorithm~\ref{algminpolyverify}.
The maximal error probability for our Monte Carlo algorithm was
$\varepsilon = 1/100$ for all runs.

\subsection{The test matrices}
Next, we describe the matrices $M_1, \ldots, M_{10}$ we used.

(a) \quad The matrices $M_1$ and $M'_1$ were purely random matrices from
$\F_3^{1000\times 1000}$ with all entries
chosen with uniform distribution from the field $\F_3$. Such matrices
are with very high probability cyclic, that is, their characteristic and
minimal polynomials are equal. Usually, Algorithm~\ref{algcharpoly} only
has to pick very few random vectors for such matrices. The {\bf for} loop
of Algorithm~\ref{algminpolyverify} quickly checks whether the least common multiple
of the relative order polynomials (which is the input candidate polynomial)
already has degree $n$. It turned out that $M_1'$ was cyclic but not $M_1$, 
and this explains the
big differences in the runtimes for these matrices.

(b)\quad The matrix $M_2$ is one coming from actual applications. Namely, it is
the matrix $a+b+ab$ where the two matrices $a,b \in \F_2^{4370 \times
4370}$ describe the action of two standard generators of the Baby monster
sporadic simple group on its smallest faithful simple module over $\F_2$.
The matrices $a$ and $b$ were downloaded from the WWW Atlas of group
representations (see \cite{WWWAtlas}). The matrix $a+b+ab$ is
interesting because it is one of the algebra words that is used
in the {\sc MeatAxe} (see \cite{MeatAxeParker} and \cite{MeatAxeHoltRees})
to compute composition series of modules and
we could very well imagine using the minimal polynomial instead of
the characteristic polynomial in some places in the {\sc MeatAxe}.

The reason why the standard algorithm for the minimal polynomial
performed rather badly on this matrix is that its characteristic polynomial
has irreducible factors of degrees $1$, $1$, $2$, $4$, $6$, $88$,
$197$, $854$ and $934$ with respective multiplicities
$2$, $2277$, $4$, $1$, $1$, $1$, $1$, $1$ and $1$.
Therefore the standard algorithm spins up large subspaces
many times.

(c)\quad The matrices $M_3$ -- $M_7$ were constructed in the following way:
In the language of $\F[x]$-modules we chose the order polynomials of
the generators of their primary cyclic submodules, that is we chose the
minimal polynomials on the primary cyclic submodules. For irreducible factors
of degree one this amounts to choosing the sizes and numbers of the
Jordan blocks occurring in the Jordan normal form of the matrix.
After writing down the corresponding normal form of the matrix we
conjugated it with a random element of the general linear group to
get a dense matrix with the same normal form.

For $M_3 \in \F_5^{600 \times 600}$ we chose one cyclic summand with
minimal polynomial $(x-\zeta_5)^{300}$ plus 300 summands with minimal
polynomial $x-\zeta_5$, where $\zeta_5 \in \F_5$ is a primitive root.
This is a typical case in which our Monte Carlo algorithm and
the deterministic verification both perform very well in comparison
with older techniques. The reason for this is that the high dimensional
cyclic subspace is spun up many times in the standard minimal polynomial
algorithm as for the matrix $M_2$.

For $M_4 \in \F_3^{1200 \times 1200}$ we chose 400 cyclic summands with
minimal polynomial $(x-\zeta_3)^2$ plus 400 cyclic summands with
minimal polynomial $(x-\zeta_3)$, where again $\zeta_3 \in \F_3$ is
a primitive root. In contrast with the matrix $M_3$,  our algorithms 
performed very well in this case but they were not
much faster than the older techniques, since the standard algorithm
run on $M_4$ does not spin up many large cyclic subspaces.

For $M_5 \in \F_{251}^{600 \times 600}$ we chose 200 different linear
factors $x-\alpha$ and for each added one cyclic space with minimal polynomial
$(x-\alpha)^2$ and one with $x-\alpha$. This example originally was a
worst case scenario for our deterministic verification. However, since $k=2$ is 
quite small, a deterministic verification can be done relatively
cheaply as described in Proposition~\ref{veryfewvectors} and
Remark~\ref{algdeterm}, even though the integer $u$ in
Algorithm~\ref{algminpolymc} is only $1$. The deterministic
verification Algorithm~\ref{algminpolyverify} ran very slowly 
(more than 300 seconds) in this example.

For $M_6 \in \F_2^{2391 \times 2391}$ we chose the irreducible polynomial
$f(x) = x^3+x^2+1 \in \F_2[x]$ of degree $3$ and added cyclic spaces with 
respective minimal polynomials $f^{400}$, $f^{200}$, $f^{100}$,
$f^{50}$, $f^{25}$, $f^{12}$, $f^6$, $f^3$ and $f$.

For $M_7 \in \F_{3^4}^{220 \times 220}$ we chose an irreducible polynomial
$f(x) \in \F_{10}[x]$ of degree $10$ and added cyclic spaces with
respective minimal polynomials $f^{10}$, $f$, $f^2$, $f^3$, $f$, $f^2$ and
$f^3$.

(d)\quad The matrices $M_8$ and $M_9$ were standard generators of 
${\rm GL}(400,17)$, 
conjugated by the pseudo-random element $M_{10}$ of the same group.
Note that $M_8$ had
order 16 while $M_9$ and $M_{10}$ had very high order and were cyclic 
matrices. We chose these examples because they may be typical of 
difficult cases in an application of the minimal polynomial algorithm
for computing the projective order of a matrix. 

Our algorithm 
very quickly discovered that the least common multiple of the relative order
polynomials was already equal to the characteristic polynomial.

\begin{figure}
\caption{Timings for minimal polynomial computation}
\label{timings}
\begin{center}
\hspace*{-5mm}
\begin{tabular}{|c|r|r|r|r|r|r|r|r|r|r|}
\hline
M & $q$ & $n$ & Lib & AS & MC & Spin & Fact & OrdP & Ver. & k \\
\hline
\hline
$M_1$  & 3   & 1000 & 1.95$^*$ & 0.65 & 13.7 & 0.33 & 13.3 & 0.05 & 0 & 2 \\
$M_1'$ & 3   & 1000 & 1.31$^*$ & 0.68 & 0.32 & 0.32 & 0 & 0 & 0 & 2 \\
$M_2$  & 2   & 4370 & 12975 & 3098 & 5.74 & 3.80 & 1.10 & 0.83 & 3.02 & 2212 \\
$M_3$  & 5   &  600 & 59.5 & 21.0 & 0.33 & 0.16 & 0.08 & 0.08 & 0.19 & 301 \\
$M_4$  & 3   & 1200 & 2.00$^*$& 0.45 & 0.44 & 0.38 & 0.06 & 0.01 & 0.06 & 800 \\
$M_5$  & 251 &  600 & 2.9 & 3.3 & 3.26 & 2.82 & 0.55 & 0.04 & 0 & 2 \\ 
$M_6$  & 2   & 2391 & 14.6 & 3.3 & 2.25  & 0.91 & 0.18 & 1.15 & 1.02 & 9 \\
$M_7$  & 243 &  220 & 0.77 & 0.88 & 0.36 & 0.34 & 0.01 & 0.01 & 0.21 & 7 \\
$M_8$  & 17  &  400 & 0.46 & 0.20  & 0.048 & 0.032 & 0.012 & 0.004 & 0.00 & 399 \\
$M_9$  & 17  &  400 & 0.26 & 0.20 & 0.23 & 0.23 & 0 & 0 & 0 & 1 \\
$M_{10}$& 17  &  400 & 0.26 & 0.19 & 0.22 & 0.22 & 0 & 0 & 0 & 1 \\
\hline
\end{tabular}
\hspace*{-5mm}

\medskip
$^*$ averaged over 10 runs
\end{center}
\end{figure}

\bigskip\noindent
{\bf Acknowledgements:}\quad 
We would like to thank an anonymous referee for invaluable
suggestions that improved and streamlined the exposition and in particular encouraged us to think about and improve
the deterministic verification procedures.

This research forms part of a Discovery
Project and Federation Fellowship Grant of the  second author funded by the
Australian Research Council.

\bibliography{minpoly}

\end{document}